\documentclass[11pt,leqno]{article}
\usepackage{amsmath, amsthm, amsfonts, amssymb}
\usepackage{mathrsfs}
\newfam\msbfam
\font\tenmsb=msbm10    \textfont\msbfam=\tenmsb \font\sevenmsb=msbm7
\scriptfont\msbfam=\sevenmsb \font\fivemsb=msbm5
\scriptscriptfont\msbfam=\fivemsb

\newfam\bigfam
\font\tenbig=msbm10 scaled \magstep2   \textfont\bigfam=\tenbig
\font\sevenbig=msbm7 scaled \magstep2 \scriptfont\bigfam=\sevenbig
\font\fivebig=msbm5 scaled \magstep2
\scriptscriptfont\bigfam=\fivebig

\def\dfrac{\displaystyle\frac}

\newtheorem{ex}{Example}[section]
\newtheorem{thm}{Theorem}[section]
\newtheorem{lem}{Lemma}[section]

\newtheorem{rem}{Remark}[section]
\newtheorem{cor}{Corollary}[section]
\newtheorem{defn}{Definition}[section]
\newtheorem{pf}{Proof}

\begin{document}

{\begin{center}\bf\LARGE Generalized $s$-Convex Functions on Fractal
Sets\end{center}}

\vspace{0.2cm}
 {\bf Huixia Mo, and Xin Sui}
 {\small

\vspace{0.2cm} School of Science, Beijing University of Posts and
Telecommunications, Beijing, 100876, China \vspace{0.2cm}

Corresponding should be address to Hui-xia Mo: huixmo@bupt.edu.cn}
\vspace{0.3cm}

\begin{minipage}{5in}

\small {We introduce two kinds of generalized $s$-convex functions
on real linear fractal sets  $\mathbb{R}^{\alpha}(0<\alpha<1)$. And
similar to the class situation, we also study the properties of
these two kinds of generalized $s$-convex functions and discuss the
relationship between them. Furthermore, some applications are
given.}
\end{minipage}

\section{Introduction}
Let $f:I\subseteq\mathbb{R}\rightarrow\mathbb{R}.$ For any $u, v\in
I$ and $t\in[0, 1],$ if the following inequality
$$f(tu+(1-t)v)\leq tf(u)+(1-t)f(v)$$
holds, then $f$ is called a convex function on $I.$

The convexity of functions plays a significant role in many fields,
such as in biological system, economy, optimization, and so on
\cite{GRA,GR}. In \cite{HL}, Hudzik and Maligranda generalized the
definition of convex function and considered, among others, two
kinds of functions which are $s$-convex.

 Let $0<s\leq 1$ and $\mathbb{R_+}=[0,\infty),$ then the two kinds of $s$-convex functions are defined respectively in the following way.

 \begin{defn} A function  $f:\mathbb{R_+}\rightarrow \mathbb{R},$ is said to be $s$-convex in the first sense if
 $$f(\alpha u+\beta v)\leq\alpha^s f(u)+\beta^s f(v),$$ for all $u,v\in\mathbb{R_+}$ and all $\alpha,\beta\geq 0$ with $\alpha^s+\beta^s=1.$ We denote this by $f\in K_s^1.$
 \end{defn}

 \begin{defn} \;A function  $f:\mathbb{R_+}\rightarrow \mathbb{R},$ is said to be $s$-convex in the second sense if
 $$f(\alpha u+\beta v)\leq\alpha^s f(u)+\beta^s f(v),$$ for all $u,v\in\mathbb{R_+}$ and all $\alpha,\beta\geq 0$ with $\alpha+\beta=1.$ We denote this by $f\in K_s^2.$
 \end{defn}

It is obvious that the $s$-convexity means just the convexity when
$s=1,$ no mater in the first sense or in the second sense. In
\cite{HL}, some properties of $s$-convex functions in both senses
are considered and various examples and counterexamples are given.
There are many research results  related to the $s$-convex
functions, see \cite{S,AK,SYQ}, and so on.

In recent years, the fractal has received significantly remarkable
attention from scientists and engineers. In the sense of Mandelbrot,
 a fractal set is the one whose Hausdorff dimension strictly exceeds the topological dimension \cite{KG,GB,M,F,E,KST}.

 The calculus on fractal set can lead to better comprehension for the various real world models from science and engineering \cite{GB}.
  Researchers have constructed many kinds of fractional calculus on fractal sets by using different approaches.
Particularly, in \cite{Y}, Yang stated the analysis of local
fractional functions on fractal space systematically, which includes
local fractional calculus and so on.   In \cite{MSY}, the authors
introduced the generalized convex function on fractal sets and
establish the generalized Jensen's inequality and generalized
Hermite-Hadamard's inequality related to generalized convex
function.  And, in \cite{WSZW}, Wei etal. established a local
fractional integral inequality on fractal space analogous to
Anderson's inequality for generalized convex functions. The
generalized convex function on fractal sets
$\mathbb{R}^{\alpha}(0<\alpha<1)$ can be stated
 as follows.

Let $f: I\subset \mathbb{R}\rightarrow\mathbb{R}^{\alpha}.$ For any
$u, v\in I$ and $t\in[0, 1],$ if the following inequality
$$f(tu+(1-t)v)\leq t^{\alpha}f(u)+(1-t)^{\alpha}f(v)$$
holds, then $f$ is called a generalized convex on $I.$

Inspired by these investigations,  we will introduce the generalized
$s$-convex function in the first or second sense on fractal sets and
study the properties of generalized $s$-convex functions.

The paper is organized as follows. In Section 2, we state the
operations with real line number fractal sets and give the
definitions of the local fractional calculus. In Section 3, we
introduce the definitions of two kinds of generalized $s$-convex
functions and study the properties of these functions.
 In Section 4, we give some applications for the two kinds of generalized $s$-convex functions on fractal sets.
\section{Preliminaries}

Let us recall the operations with real line number on fractal space
and use the Gao-Yang- Kang's idea to describe the definitions of the
 local fractional derivative and local fractional integral \cite{Y,YBM,GYK,YG1,YBKD}.

If $a^\alpha,b^\alpha$ and $c^\alpha$ belong to the set
$\mathbb{R}^\alpha(0<\alpha\leq1)$ of real line numbers, then one
has the following:

\vspace{0.2cm}
 \noindent
(1)\;\;$a^\alpha+b^\alpha$ and $a^\alpha b^\alpha$ belong to the set $\mathbb{R}^\alpha$;\\
(2)\;\;$a^\alpha+b^\alpha=b^\alpha+a^\alpha=(a+b)^\alpha=(b+a)^\alpha$;\\
(3)\;\;$a^\alpha+(b^\alpha+c^\alpha)=(a^\alpha+b^\alpha)+c^\alpha$;\\
(4)\;\;$a^\alpha b^\alpha=b^\alpha a^\alpha=(ab)^\alpha=(ba)^\alpha$;\\
(5)\;\;$a^\alpha(b^\alpha c^\alpha)=(a^\alpha b^\alpha)c^\alpha$;\\
(6)\;\;$a^\alpha(b^\alpha+c^\alpha)=a^\alpha b^\alpha+a^\alpha c^\alpha$;\\
(7)\;\;$a^\alpha+0^\alpha=0^\alpha+a^\alpha=a^\alpha$ and
$a^\alpha\cdot1^\alpha=1^\alpha\cdot a^\alpha=a^\alpha$.

Let us now state some definitions about the local fractional
calculus on $\mathbb{R}^\alpha.$

\begin{defn}\cite{Y}\;\; A non-differentiable function
$f: \mathbb{R} \rightarrow \mathbb{R}^\alpha,\;x \to f(x)$ is called
to be local fractional continuous at $x_0,$ if for any
$\varepsilon>0,$ there exists $\delta>0,$ such that
$$ |f(x)-f(x_0)|<\varepsilon^\alpha$$
holds for $|x-x_0|<\delta,$ where $\varepsilon,
\delta\in\mathbb{R}.$ If $f$ is local fractional continuous on the
interval $(a,b),$ one denotes $f\in C_{\alpha}(a,b).$
\end{defn}

\begin{defn}\cite{Y} The local fractional derivative of function $f$ of order $\alpha$ at $x=x_0$ is defined by
$$f^{(\alpha)}(x_0)=\frac{d^\alpha f(x)}{dx^\alpha}|_{x=x_0}=\lim\limits_{x\rightarrow x_0}\frac{\triangle^\alpha(f(x)-f(x_0))}{(x-x_0)^\alpha},$$
where $\triangle^\alpha(f(x)-f(x_0))=\Gamma(1+a)(f(x)-f(x_0))$ and
the Gamma function is defined by
$\Gamma(t)=\int_{0}^{+\infty}x^{t-1}e^{-x}dx.$

If there exists $f^{((k+1)\alpha)}(x)=\overbrace{D_x^{\alpha} \ldots
D_x^\alpha}\limits^{k+1\; \mbox {times}}f(x)$ for any $x\in
I\subseteq\mathbb{R},$  then we denoted $f\in D_{(k+1)\alpha}(I)$,
where $k=0,1,2\dots.$
\end{defn}

\begin{defn}\cite{Y}\;\;Let $f\in C_{\alpha}[a,b]$. Then the local fractional integral of the function $f$  of order $\alpha$ is defined by,
$$\begin{array}{cl}&_aI_b^{(\alpha)}f\\
& =\dfrac{1}{\Gamma(1+a)}\int_a^bf(t)(dt)^\alpha\\&
=\dfrac{1}{\Gamma(1+a)}\lim\limits_{\triangle t\rightarrow
0}\sum\limits_{j=0}^{N}f(t_j)(\bigtriangleup
t_j)^\alpha,\end{array}$$ with $\triangle t_j=t_{j+1}-t_j,$
$\triangle t=\max\{\triangle t_1,\triangle t_2,\triangle
t_j,\cdots,\triangle t_N-1\}$ and $[t_j,t_j+1],j=0,\cdots,N-1,$
where $t_0=a<t_1<\cdots<t_i<\cdots<t_N=b$ is a partition of the
interval $[a,b]$.
\end{defn}

\begin{lem}\cite{Y}\;\;Suppose that $f, g\in C_{\alpha}[a,b]$ and
$f, g\in D_{\alpha}(a,b)$. If $\lim\limits_{x\rightarrow
x_0}f(x)=0^\alpha,$ $\lim\limits_{x\rightarrow x_0}g(x)=0^\alpha$
and $g^{(\alpha)}(x)\neq0^\alpha.$ Suppose that
$\lim\limits_{x\rightarrow
x_0}\dfrac{f^{(\alpha)}(x)}{g^{(\alpha)}(x)}=A^\alpha,$ then
$$\lim\limits_{x\rightarrow x_0}\dfrac{f(x)}{g(x)}=A^\alpha.$$
\end{lem}

\begin{lem}\cite{Y}\;\;Suppose that $f(x)\in C_{\alpha}[a,b]$, then

$$\dfrac{d^{\alpha}(_aI_x^{(\alpha)}f)}{dx^{\alpha}}=f(x),\;\;a<x<b.$$
\end{lem}

\section{Generalized s-convexity functions}

The convexity of functions plays a significant role in many fields.
In this section, let us introduce  two kinds of  generalized
s-convex functions on fractal sets. And then, we study the
properties of the two kinds of  generalized s-convex functions.

\begin{defn}Let $\mathbb{R_+}=[0, +\infty).$ A function $f:\mathbb{R_+}\rightarrow\mathbb{R^{\alpha}}$
is said to be generalized $s$-convex $(0<s<1)$ in the first sense, if
$$f(\lambda_1u+\lambda_2v)\leq\lambda_1^{s\alpha}f(u)+\lambda_2^{s\alpha}f(v), \eqno (3.1)$$
for all $u,v\in\mathbb{R_+}$ and all $\lambda_1,\lambda_2\geq0$ with $\lambda_1^ s+\lambda_2^s=1.$ We denote this by $f\in GK_s^1.$\end{defn}

\begin{defn}A function $f:\mathbb{R_+}\rightarrow\mathbb{R^{\alpha}}$
is said to be generalized $s$-convex $(0<s<1)$ in the second sense, if
$$f(\lambda_1u+\lambda_2v)\leq\lambda_1^{s\alpha}f(u)+\lambda_2^{s\alpha}f(v), \eqno (3.2)$$
for all $u,v\in\mathbb{R_+}$ and all $\lambda_1,\lambda_2\geq0$ with $\lambda_1 +\lambda_2=1.$ We denote this by $f\in GK_s^2.$
\end{defn}

Note that, when $s=1$, the generalized $s$-convex functions in both
sense are the generalized  convex functions, see \cite{MSY}.

\begin{thm}
Let $0<s<1.$\\
(a)\;\;If $f\in GK_s^1,$ then $f$ is non-decreasing on $(0,
+\infty)$ and $$f(0^+)=\lim\limits_{u\rightarrow 0^+}f(u)\leq
f(0).$$ (b)\;\;If $f\in GK_s^2,$ then $f$ is non-negative on
$[0,+\infty).$
\end{thm}

\begin{pf}
(a)\;\;Since $f\in GK_s^1,$ we have, for $u>0$ and $\lambda\in[0,1],$
$$f[(\lambda^{1/s}+(1-\lambda)^{1/s})u]\leq\lambda^\alpha f(u)+(1-\lambda)^\alpha f(u)=f(u).$$
The function $h(\lambda)=\lambda^{1/s}+(1-\lambda)^{1/s},$ is continuous on $[0,1],$  decreasing on $[0,1/2],$
 increasing on $[1/2,1]$ and $h([0,1])=[h(1/2),h(1)]=[2^{1-1/s},1].$ This yields that $$f(tu)\leq f(u),\eqno (3.3)$$ for $u>0$ and $t\in[2^{1-1/s},1].$
  If $t\in[2^{2(1-1/s)},1],$ then $t^{1/2}\in[2^{1-1/s},1].$ Therefore, by the fact that (3.3) holds for all $u>0,$ we get
$$f(tu)=f(t^{1/2}(t^{1/2}u))\leq f(t^{1/2}u)\leq f(u),$$ for all $u>0.$ So we can obtain that $$f(tu)\leq f(u), \;\mbox{for all}\; u>0, t\in (0,1].$$

So, taking $0<u<v,$ we get $$f(u)=f((u/v)v)\leq f(v),$$ which means
that $f$ is non-decreasing on $(0, +\infty).$

As for the second part, for $u>0$ and $\lambda_1, \lambda_2\geq0$
 with $\lambda_1^ s+\lambda_2^s=1,$ we have
$$f(\lambda_1u)=f(\lambda_1u+\lambda_20)\leq\lambda_1^{s\alpha}f(u)+\lambda_2^{s\alpha}f(0).$$
And taking $u\rightarrow0^+,$ we get $$\lim\limits_{u\rightarrow
0^+}f(u)=\lim\limits_{u\rightarrow
0^+}f(\lambda_1u)\leq\lambda_1^{s\alpha}\lim\limits_{u\rightarrow
0^+}f(u)+\lambda_2^{s\alpha}f(0).$$ So, $$\lim\limits_{u\rightarrow
0^+}f(u)\leq f(0).$$ (b)\;\; For $f\in GK_s^2,$ we can get that for
$u\in \mathbb{R_+},$ $$f(u)=f(u/2+u/2)\leq f(u)/2^{s\alpha
}+f(u)/2^{s\alpha}=2^{(1-s)\alpha}f(u).$$

So, $(2^{1-s}-1)^\alpha f(u)\geq0^\alpha.$ This means that $f(u)\geq0^\alpha,$ since $0<s<1.$
\end{pf}

\begin{rem}The above results do not hold, in general, in the case of generalized convex functions, i.e. when $s=1$,
 because a  generalized convex function $f:\mathbb{R_+}\rightarrow \mathbb{R^{\alpha}},$ need not be either non-decreasing or non-negative.
\end{rem}

\begin{rem}
If  $0<s<1,$ then the function $f\in GK_s^1$ is non-decreasing on
$(0, +\infty)$ but not necessarily on $[0, +\infty).$
\end{rem}
Function $F:\mathbb{R}^2\rightarrow \mathbb{R}^\alpha$ is called to
be generalized convex in each variable, if
$$ F(\lambda_1^{s}u+\lambda_2^{s}v, \lambda_1^{s}r+\lambda_2^{s}t)\leq\lambda_1^{s\alpha}F(u,r)+\lambda_2^{s\alpha}F(v,t).$$
 for all $(u,v),$  $(r,t)\in \mathbb{R}^2$ and $\lambda_1,\lambda_2\in[0,1].$

\begin{thm}\;\;
Let $0<s<1.$ If $f,g:\mathbb{R}\rightarrow\mathbb{R}$ and $f, g\in
K_s^1$ and if $F:\mathbb{R}^2\rightarrow \mathbb{R}^\alpha$ is a
generalized convex and non-decreasing function in each variable,
then the function $h:\mathbb{R_+}\rightarrow \mathbb{R}^\alpha$
defined by
$$h(u)=F(f(u),g(u))$$ is in $GK_s^1.$ In particular, if $f, g\in
K_s^1$ then $f^{\alpha}+g^{\alpha},$ $\max\{f^{\alpha},
g^{\alpha}\}\in GK_s^1.$
\end{thm}
\begin{pf}If $u,v\in \mathbb{R_+},$ then for all $\lambda_1,\lambda_2\geq 0$ with $\lambda_1^s+\lambda_2^s=1,$
$$\begin{array}{cl}h(\lambda_1u+\lambda_2v)
& =F(f(\lambda_1u+\lambda_2v),g(\lambda_1u+\lambda_2v))\\& \leq
F(\lambda_1^{s}f(u)+\lambda_2^{s}f(v),
\lambda_1^{s}g(u)+\lambda_2^{s}g(v))\\& \leq\lambda_1^{s
\alpha}F(f(u),g(u))+\lambda_2^{s\alpha}F(f(v),g(v))\\
& =\lambda_1^{s\alpha}h(u)+\lambda_2^{s\alpha}h(v).\end{array}$$
Thus, $h\in GK_s^{1}.$

Moreover, since $F(u,v)=u^\alpha+v^\alpha,$ $F(u,v)=\max\{u^\alpha,
v^\alpha\}$ are non-decreasing generalized convex functions on
$R^2,$ so they yield particular cases of our theorem.\end{pf}

Let's pay attention to the situation when the condition
$\lambda_1^s+\lambda_2^s=1$ $(\lambda_1+\lambda_2=1)$ in the
definition of $GK_s^1(GK_s^2)$ can be equivalently replaced by the
condition $\lambda_1^s+\lambda_2^s\leq 1$ $(\lambda_1+\lambda_2\leq
1).$

\begin{thm}
(a)\;\;Let $f\in GK_s^1.$ Then inequality (3.1) holds for all $u,v\in R_+$ and all $\lambda_1,\lambda_2\geq 0$ with $\lambda_1^s+\lambda_2^s< 1$ if and only if  $f(0)\leq 0^\alpha.$\\
(b)\;\;Let $f\in GK_s^2.$ Then inequality (3.2) holds for all $u,v\in R_+$ and all $\lambda_1,\lambda_2\geq0$ with $\lambda_1+\lambda_2< 1$ if and only if $f(0)= 0^\alpha.$
\end{thm}
\begin{pf}
(a)\;\;Necessity is obvious by taking $u=v=0$ and $\lambda_1=\lambda_2=0.$ Let us show the sufficiency.

Assume that $u,v\in \mathbb{R_+}$ and $\lambda_1,\lambda_2\geq0$
with  $0<\lambda_3=\lambda_1^s+\lambda_2^s<1.$ Put
$a=\lambda_1\lambda_3^{-1/s}$ and $b=\lambda_2\lambda_3^{-1/s}.$
Then $a^s+b^s=\lambda_1^s/\lambda_3+\lambda_2^s/\lambda_3=1$ and
$$\begin{array}{cl}f(\lambda_1u+\lambda_2v)
& =f(a\lambda_3^{1/s}u+b\lambda_3^{1/s}v)\\& \leq a^{s
\alpha}f(\lambda_3^{1/s}u)+b^{s \alpha}f(\lambda_3^{1/s}v)\\& =a^{s
\alpha}f[\lambda_3^{1/s}u+(1-\lambda_3)^{1/s}0]+b^{s \alpha
}f[\lambda_3^{1/s}v+(1-\lambda_3)^{1/s}0]\\& \leq a^{
 s\alpha}[\lambda_3^\alpha f(u)+(1-\lambda_3)^\alpha f(0)]+b^{
 s\alpha}[\lambda_3^\alpha f(v)+(1-\lambda_3)^\alpha f(0)]\\& =a^{
 s\alpha }\lambda_3^\alpha f(u)+b^{s \alpha}\lambda_3^\alpha
f(v)+(1-\lambda_3)^\alpha f(0)\\& \leq \lambda_1^{s\alpha
}f(u)+\lambda_2^{s \alpha}f(v).\end{array}$$ (b)\;\;Necessity.

 Taking $u=v=\lambda_1=\lambda_2=0,$ we obtain $f(0)\leq 0^\alpha.$ And using Theorem 3.1(b), we get $f(0)\geq0^\alpha.$ Therefore $f(0)=0^\alpha.$

Sufficiency.

Let $u,v\in \mathbb{R_+}$ and $\lambda_1,\lambda_2\geq0$ with
$0<\lambda_3=\lambda_1+\lambda_2<1.$ Put $a=\lambda_1/\lambda_3$ and
$b=\lambda_2/\lambda_3,$ then $a+b=1.$

 So,
$$\begin{array}{cl}f(\lambda_1u+\lambda_2v)
& =f(a\lambda_3u+b\lambda_3v)\\
& \leq a^{s \alpha}f(\lambda_3u)+b^{s
\alpha}f(\lambda_3v)\\
& =a^{s\alpha}f[\lambda_3u+(1-\lambda_3)0]+b^{s\alpha
}f[\lambda_3v+(1-\lambda_3)0]\\& \leq a^{s \alpha}[\lambda_3^{s
\alpha }f(u)+(1-\lambda_3)^{s \alpha}f(0)]+b^{s \alpha}[\lambda_3^{s
\alpha }f(v)+(1-\lambda_3)^{s \alpha}f(0)]\\& =a^{s\alpha
}\lambda_3^{s \alpha}f(u)+b^{s \alpha}\lambda_3^{s\alpha
}f(v)+(1-\lambda_3)^{s
\alpha}f(0)\\
&= \lambda_1^{s \alpha}f(u)+\lambda_2^{s\alpha}f(v).\end{array}$$
\end{pf}

\begin{thm}
(a)\;\;Let $0<s\leq1.$ If $f\in GK_s^2$ and $f(0)=0^\alpha,$ then $f\in GK_s^1.$\\
(b)\;\;Let $0<s_1\leq s_2\leq 1.$ If $f\in GK_{s_2}^2$ and $f(0)=0^\alpha,$ then $f\in GK_{s_1}^2.$\\
(c)\;\;Let $0<s_1\leq s_2\leq 1.$ If $f\in GK_{s_2}^1$ and $f(0)\leq 0^\alpha,$ then  $f\in GK_{s_1}^1.$
\end{thm}

\begin{pf}

(a)\;\;Assume that $f\in GK_s^2$ and $f(0)=0^\alpha.$ Let $\lambda_1,\lambda_2\geq 0$ with $\lambda_1^s+\lambda_2^s=1,$
 we have $\lambda_1+\lambda_2\leq \lambda_1^s+\lambda_2^s=1.$ From Theorem 3.3(b), we can get
$$f(\lambda_1u+\lambda_2v)\leq \lambda_1^{s \alpha}f(u)+\lambda_2^{s \alpha}f(v),$$ for $u,v\in \mathbb{R_+},$ then $f\in GK_s^1.$

(b)\;\;Assume that $f\in GK_{s_2}^2,$ $u,v\in \mathbb{R_+}$ and
$\lambda_1,\lambda_2\geq 0$ with $\lambda_1+\lambda_2=1.$ Then we
have
$$f(\lambda_1u+\lambda_2v)\leq\lambda_1^{s_2
\alpha}f(u)+\lambda_2^{s_2 \alpha}f(v)\leq\lambda_1^{s_1\alpha
}f(u)+\lambda_2^{s_1 \alpha}f(v).$$
So $f\in GK_{s_1}^2.$

(c)\;\;Assume that $f\in GK_{s_2}^1,$ $u,v\in R_+$ and
$\lambda_1,\lambda_2\geq 0$ with
$\lambda_1^{s_1}+\lambda_2^{s_1}=1.$ Then
$\lambda_1^{s_2}+\lambda_2^{s_2}\leq
\lambda_1^{s_1}+\lambda_2^{s_1}=1.$ Thus, according to Theorem
3.3(a), we have
$$f(\lambda_1u+\lambda_2v)\leq\lambda_1^{s_2\alpha}f(u)+\lambda_2^{s_2\alpha}f(v)\leq\lambda_1^{s_1\alpha}f(u)+\lambda_2^{s_1 \alpha}f(v).$$
So, $f\in GK_{s_1}^1.$
\end{pf}

\begin{thm}Let $0<s<1$ and  $p:\mathbb{R_+}\rightarrow \mathbb{R^\alpha_+}$ be a non-decreasing function. Then the function $f$ defined for $u\in \mathbb{R_+}$ by
$$f(u)=u^{(s/(1-s))\alpha }p(u)$$
belongs to $GK_s^1.$
\end{thm}
\begin{pf}Let $v\geq u\geq 0$ and $\lambda_1,\lambda_2\geq 0$ with $\lambda_1^s+\lambda_2^s=1.$ We consider two cases.

{\bf Case \uppercase\expandafter{\romannumeral1}}:\;\; It is easy to
see that $f$ is a non-decreasing function. Let
$\lambda_1u+\lambda_2v\leq u,$ then
$$f(\lambda_1u+\lambda_2v)\leq f(u)=(\lambda_1^{s \alpha}+\lambda_2^{s \alpha})f(u)\leq\lambda_1^{s \alpha}f(u)+\lambda_2^{s \alpha}f(v).$$

{\bf Case \uppercase\expandafter{\romannumeral2}}:\;\;Let
$\lambda_1u+\lambda_2v>u,$ then $\lambda_2v>(1-\lambda_1)u.$ So,
$\lambda_2>0$ and $\lambda_1\leq\lambda_1^s.$ Thus,
$$\lambda_1-\lambda_1^{s+1}\leq\lambda_1^s-\lambda_1^{s+1},$$ i.e.,
$$\lambda_1/(1-\lambda_1)\leq\lambda_1^s/(1-\lambda_1^s)=(1-\lambda_2^s)/\lambda_2^s,$$ and
$$\lambda_1\lambda_2/(1-\lambda_1)\leq\lambda_2^{1-s}-\lambda_2.$$

Thus, we can get that
$$\lambda_1u+\lambda_2v\leq(\lambda_1+\lambda_2)v\leq(\lambda_1^s+\lambda_2^s)v=v,$$
and
$$\lambda_1u+\lambda_2v\leq\lambda_1\lambda_2v/(1-\lambda_1)+\lambda_2v\leq(\lambda_2^{1-s}-\lambda_2)v+\lambda_2v=\lambda_2^{1-s}v.$$
Then, $$(\lambda_1u+\lambda_2v)^{s/(1-s)}\leq\lambda_2^sv^{s/(1-s)}.$$

We obtain
$$\begin{array}{cl}f(\lambda_1u+\lambda_2v)
&=(\lambda_1u+\lambda_2
v)^{(s/(1-s))\alpha}p(\lambda_1u+\lambda_2v)\leq\lambda_2^{s
\alpha}v^{(s/(1-s))\alpha}p(v)\\
&=\lambda_2^{s\alpha}f(v)\leq\lambda_1^{s
\alpha}f(u)+\lambda_2^{s\alpha}f(v).\end{array}$$
\end{pf}

\begin{thm}
(a)\;\;Let $f\in GK_{s_1}^1$ and $g\in K_{s_2}^1,$ where $0<s_1,s_2\leq 1.$ If $f$ is a non-decreasing function and $g$ is a non-negative function such that $f(0)\leq 0^\alpha$ and  $g(0)=0,$ then the composition $f\circ g$ of $f$ with $g$ belongs to $GK_s^1,$ where $s=s_1s_2.$\\
(b)\;\;Let $f\in GK_{s_1}^1$ and $g\in GK_{s_2}^1,$ where
$0<s_1,s_2\leq 1.$ Assume that $0<s_1,s_2<1.$ If $f$ and $g$ are
non-negative functions such that either $f(0)=0^\alpha$ and
 $g(0^+)=g(0),$ or $g(0)=0^\alpha$ and $f(0^+)=f(0),$ then the product $fg$ of $f$ and $g$ belongs to $GK_s^1,$ where $s=\min\{s_1,s_2\}.$
\end{thm}
\begin{pf}
(a)\;\;Let $u,v\in \mathbb{R_+},$ $\lambda_1,\lambda_2\geq0$ with
$\lambda_1^s+\lambda_2^s=1,$ where $s=s_1s_2.$ Since
$\lambda_1^{s_i}+\lambda_2^{s_i}\leq\lambda_1^{s_1s_2}+\lambda_2^{s_1s_2}=1$
for $i=1,2,$ then according to  Theorem 3 (a) in \cite{HL} and
Theorem 3.3(a) in the paper, we have
$$\begin{array}{cl}f\circ g(\lambda_1u+\lambda_2v)&=f(g(\lambda_1u+\lambda_2v))\\
& \leq f(\lambda_1^{s_2}g(u)+\lambda_2^{s_2}g(v))\\&
\leq\lambda_1^{s_1s_2\alpha}f(g(u))+\lambda_2^{
s_1s_2\alpha}f(g(v))\\
& =\lambda_1^{s\alpha}f\circ g(u)+\lambda_2^{s\alpha}f\circ
g(v),\end{array}$$
which means that $f\circ g\in GK_s^1.$\\
(b)\;\;According to Theorem 3.1(a), $f,g$ are non-decreasing on $(0,
+\infty).$

So,
$$(f(u)-f(v))(g(v)-g(u))\leq 0^{\alpha},$$
or, equivalently,
$$f(u)g(v)+f(v)g(u)\leq f(u)g(u)+f(v)g(v),$$
for all $v>u>0.$

If $v>u=0,$ then the inequality is still true because $f,g$ are non-negative and either $f(0)=0^\alpha$ and $g(0^+)=g(0)$ or $g(0)=0^\alpha$ and $f(0^+)=f(0).$

Now let $u,v\in \mathbb{R_+}$ and $\lambda_1,\lambda_2\geq 0$ with
$\lambda_1^s+\lambda_2^s=1,$ where $s=\min\{s_1,s_2\}.$ Then
$\lambda_1^{s_i}+\lambda_2^{s_i}\leq\lambda_1^s+\lambda_2^s=1$ for
$i=1,2.$ And by Theorem 3.3(a), we have
$$\begin{array}{cl}&f(\lambda_1u+\lambda_2v)g(\lambda_1u+\lambda_2v)\\
&\leq(\lambda_1^{s_1\alpha}f(u)+\lambda_2^{s_1\alpha}f(v))(\lambda_1^{s_2 \alpha}g(u)+\lambda_2^{s_2 \alpha }g(v))\\
&=\lambda_1^{(s_1+s_2)\alpha}f(u)g(u)+\lambda_1^{s_1 \alpha
}\lambda_2^{s_2 \alpha}f(u)g(v)+\lambda_1^{ s_2\alpha}\lambda_2^{s_1 \alpha}f(v)g(u)+\lambda_2^{(s_1+s_2)\alpha}f(v)g(v)\\
&\leq\lambda_1^{2s\alpha }f(u)g(u)+\lambda_1^{s \alpha}\lambda_2^{s
\alpha}(f(u)g(v)+f(v)g(u))+\lambda_2^{2s \alpha}f(v)g(v)\\
&\leq\lambda_1^{2s \alpha}f(u)g(u)+\lambda_1^{s \alpha}\lambda_2^{s
\alpha}(f(u)g(u)+f(v)g(v))+\lambda_2^{2s \alpha}f(v)g(v)\\
&=\lambda_1^{s \alpha}f(u)g(u)+\lambda_2^{s \alpha}f(v)g(v),\end{array}$$\\
which means that $fg\in GK_s^1.$
\end{pf}

\begin{rem}From the above proof, we can get that, if $f$ is a non-decreasing function in $GK_s^2$ and $g$ is
a non-negative convex function on $[0, +\infty),$ then the
composition $f\circ g$ of $f$ with $g$ belongs to $GK_s^2.$
\end{rem}

\begin{rem}
Generalized convex functions on $[0, +\infty)$ need not be
monotonic. However, if $f$ and $g$ are non-negative, generalized
convex and either both are non-decreasing or both are non-increasing
on $[0, +\infty)$ then the product $fg$ is also a generalized convex
function.\end{rem}

Let  $f:\mathbb{R_+}\rightarrow\mathbb{R}_+$ be a continuous
function. Then $f$ is said to be a $\varphi$-function if $f(0)=0$
and  $f$ is non-decreasing on $\mathbb{R_+}.$ Similarly, we can
define the  $\varphi$-type function on fractal sets as follows: A
function $f:\mathbb{R_+}\rightarrow\mathbb{R}_+^\alpha$ is said to
be a  $\varphi$-type function if $f(0)=0^\alpha$ and $f\in
C_{\alpha}(\mathbb{R_+})$ is non-decreasing.

\begin{cor}If $\Phi$ is a generalized convex $\varphi$-type function and $g\in K_s^1$ is a $\varphi$-function, then the composition $\Phi \circ g$ belongs to $GK_s^1.$ In
particular, the $\varphi$-type function $h(u)=\Phi(u^s)$ belongs to
$GK_s^1.$
\end{cor}

\begin{cor}If $\Phi$ is a convex  $\varphi$-function and $f\in GK_s^2$ is a $\varphi$-type function, then
 the composition $f\circ \Phi$ belongs to $GK_s^2.$ In particular, the $\varphi$-type function $h(u)=[\Phi(u)]^{s\alpha}$ belongs to $GK_s^2.$
\end{cor}
\begin{thm}If $0<s<1$ and $f\in GK_s^1$ is a  $\varphi$-type function, then there exits a generalized convex $\varphi$-type function $\Phi$ such that $$f(2^{-1/s}u)\leq\Phi(u^s)\leq f(u),$$ for all $u\geq0.$\end{thm}
\begin{pf} By the generalized $s$-convexity of the function $f$ and by $f(0)=0^\alpha,$ we obtain $f(\lambda_1u)\leq\lambda_1^{s \alpha}f(u)$ for all $u\geq0$ and all $\lambda_1\in[0,1].$

Assume now that $v>u>0$. Then $$f(u^{1/s})\leq
f((u/v)^{1/s}v^{1/s})\leq(u^\alpha/v^\alpha)f(v^{1/s}),$$
i.e.,$$f(u^{1/s})/u^\alpha\leq f(v^{1/s})/v^\alpha.\eqno (3.4)$$

Inequality (3.4) means that the function $f(u^{1/s})/u^\alpha$ is
a non-decreasing function on $(0, \infty).$ And, since $f$ is a
 $\varphi$-type function, thus $f$ is local fractional
continuous $[0,+\infty).$

 Define
$$\Phi(u)=\left\{\begin{aligned}&0^\alpha,\;\;\;\;\;\;\;\;\;\;\;\;\;\;\;\;\;\;\;\;\;\;\;\;\;\;\;\;\;\;\;\;\;\;\;\;u=0,\\&\Gamma(1+\alpha) _0I_u^{(\alpha)}(f(t^{1/s})/t^\alpha), u>0.\end{aligned}\right.$$

From Lemma 2.1 and Lemma 2.2, it is easy to see that $\Phi$ is a
generalized convex $\varphi$-type function and
 $$\Phi(u^s)=\Gamma(1+\alpha)_0I_{u^s}^{(\alpha)}(f(t^{1/s})/t^\alpha)\leq (f((u^s)^{1/s})/u^{s \alpha})u^{s \alpha}=f(u).$$
 Moreover,
 $$\Phi(u^s)\geq\Gamma(1+\alpha)_{(u^s/2)}I_{u^s}^{(\alpha)}(f(t^{1/s})/t^\alpha)\geq(f((u^s/2)^{1/s})2^\alpha u^{-s \alpha})u^{ s \alpha}/2^\alpha=f(2^{-1/s}u).$$

 Therefore,
 $$f(2^{-1/s}u)\leq\Phi(u^s)\leq f(u),$$
 for all $u\geq0.$
\end{pf}

\section{Applications}
Based on the properties of the two kinds of generalized s-convex
functions in the above section, some applications are given.

\begin{ex}Let $0<s<1,$ and $a^\alpha,b^\alpha,c^\alpha\in \mathbb{R^\alpha}.$ For $u\in \mathbb{R_+},$ define
$$ f(u)=\left\{\begin{aligned}&a^\alpha,\;\;\;\;\;\;\;\;\;\;\;\;\;u=0,\\&b^\alpha u^{s \alpha}+c^\alpha , u>0.\end{aligned}\right.$$
We have the following conclusions.\\
(i)\;\;\;If $b^\alpha\geq0^\alpha$ and $c^\alpha\leq a^\alpha,$ then $f\in GK_s^1.$\\
(ii)\;\;If $b^\alpha\geq0^\alpha$ and $c^\alpha< a^\alpha,$ then $f$ is non-decreasing on $(0, +\infty)$ but not on $[0, +\infty).$\\
(iii)\;If $b^\alpha\geq0^\alpha$ and $0^\alpha\leq c^\alpha\leq a^\alpha,$ then $f\in GK_s^2.$\\
(iv)\;If $b^\alpha>0^\alpha$ and $c^\alpha<0^\alpha,$ then $f\not\in GK_s^2.$
\end{ex}

\begin{pf}
(i)\;\; Let $\lambda_1, \lambda_2\geq0$ with
$\lambda_1^s+\lambda_2^s=1.$ Then, there are two non-trivial cases.

{\bf Case I}\;\; Let $u,v>0.$ Then $\lambda_1u+\lambda_2v>0.$

Thus,
$$\begin{array}{cl}f(\lambda_1u+\lambda_2v)&

=b^\alpha(\lambda_1u+\lambda_2v)^{s \alpha}+c^\alpha\leq
b^\alpha(\lambda_1^{s \alpha}u^{s \alpha}+\lambda_2^{s \alpha}v^{s
\alpha})+c^\alpha\\& =b^\alpha(\lambda_1^{s \alpha}u^{s \alpha
}+\lambda_2^{s \alpha}v^{s \alpha})+c^\alpha(\lambda_1^{s \alpha
}+\lambda_2^{s \alpha})\\& =\lambda_1^{s \alpha}(b^\alpha u^{s
\alpha }+c^\alpha)+\lambda_2^{s \alpha}(b^\alpha v^{s
\alpha}+c^\alpha)\\& =\lambda_1^{s \alpha}f(u)+\lambda_2^{s
\alpha}f(v).\end{array}$$

{\bf Case II}\;\; Let $v>u=0.$ We need only to consider
$\lambda_2>0.$

Thus, we have
$$\begin{array}{cl}f(\lambda_10+\lambda_2v)&
=f(\lambda_2v)=b^\alpha\lambda_2^{s \alpha}v^{s \alpha
}+c^\alpha=b^\alpha\lambda_2^{s \alpha}v^{s \alpha
}+c^\alpha(\lambda_1^{s \alpha}+\lambda_2^{s \alpha})\\&
=\lambda_1^{s \alpha}c^\alpha+\lambda_2^{s \alpha}(b^\alpha v^{s
\alpha}+c^\alpha)\\& =\lambda_1^{s \alpha}c^\alpha+\lambda_2^{s
\alpha }f(v)\leq\lambda_1^{s \alpha}a^\alpha+\lambda_2^{s
\alpha}f(v)\\& =\lambda_1^{s \alpha}f(0)+\lambda_2^{s
\alpha}f(v).\end{array}$$ So, $f\in GK_s^1.$

(ii)\;\;From Theorem 3.1, we can see that  property (ii) is true.

(iii)\;\;Let $\lambda_1, \lambda_2\geq0$ with
$\lambda_1+\lambda_2=1.$ Similar to the estimate of (i), there are
also two cases.

 Let $v,v>0.$ Then $\lambda_1u+\lambda_2v>0,$

 Thus,
$$\begin{array}{cl}f(\lambda_1u+\lambda_2v)
& =b^\alpha(\lambda_1u+\lambda_2v)^{s \alpha}+c^\alpha \\
&< b^\alpha(\lambda_1^{s \alpha}u^{s \alpha}+\lambda_2^{s
\alpha}v^{s
\alpha})+c^\alpha(\lambda_1^\alpha+\lambda_2^\alpha)\\
&\leq
b^\alpha(\lambda_1^{s \alpha}u^{s \alpha}+\lambda_2^{s \alpha}v^{s
\alpha})+c^\alpha(\lambda_1^{s \alpha}+\lambda_2^{s \alpha})
\\& =\lambda_1^{s \alpha}(b^\alpha u^{s \alpha
}+c^\alpha)+\lambda_2^{s \alpha}(b^\alpha v^{s \alpha}+c^\alpha)\\&
\leq\lambda_1^{s \alpha} f(u)+\lambda_2^{s \alpha}
f(v).\end{array}$$

 Let $v>u=0.$ We need only to consider
$\lambda_2>0.$

Thus, we have
$$\begin{array}{cl}f(\lambda_10+\lambda_2v)
& =f(\lambda_2v)=b^\alpha\lambda_2^{s \alpha}v^{s\alpha
}+c^\alpha(\lambda_1^\alpha+\lambda_2^\alpha)\\
&< b^\alpha\lambda_2^{s \alpha}v^{s\alpha}+c^\alpha(\lambda_1^{s
\alpha}+\lambda_2^{s
\alpha})\\
& =\lambda_2^{s \alpha}(b^\alpha v^{s \alpha
}+c^\alpha)+c^\alpha\lambda_1^{s \alpha}\\& \leq\lambda_2^{s
\alpha}(b^\alpha v^{s \alpha}+c^\alpha)+a^\alpha\lambda_1^{s
\alpha}\\& =\lambda_2^{s \alpha} f(v)+\lambda_1^{s \alpha}
f(0).\end{array}$$

So, $f\in GK_s^2.$

(iv)\;\;Assume that $f\in GK_s^2,$  then $f$ is non-negative on
$(0,\infty).$  On the other hand, we can take $u_1>0,c_1<0$ such
that $f(u_1)=b^\alpha u_1^{s \alpha}+c_1^\alpha<0^{\alpha},$  which
contradict the assumption.
\end{pf}

\begin{ex}Let $0<s<1$ and $k>1.$  For $u\in R_+,$ define,
$$f(u)=\left\{\begin{aligned}&u^{(s/(1-s))\alpha},\;\;\;\;\;\;0\leq u\leq 1,\\&k^\alpha u^{(s/(1-s))\alpha}, \;\; u>1.\end{aligned}\right.$$
The function $f$ is non-negative, not local fractional continuous at
$u=1$ and belongs to $GK_s^1$ but not to $GK_s^2.$\end{ex}

\begin{pf} From Theorem 3.5, we have that $f\in GK_s^1.$ In the following, Let us show that $f$ isn't in $GK_s^2.$

Take an arbitrary $a>1$ and put $u=1.$ Consider all $v>1$ such that $\lambda_1u+\lambda_2v=\lambda_1+\lambda_2v=a,$ where $\lambda_1,\lambda_2\geq 0$ and $\lambda_1+\lambda_2=1.$

If $f\in GK_s^2,$ it must be
$$k^\alpha a^{(s/(1-s))\alpha}\leq\lambda_1^{s \alpha}+k^\alpha(1-\lambda_1)^{s \alpha}[(a-\lambda_1)/(1-\lambda_1)]^{(s/(1-s))\alpha},\eqno (3.5)$$
for all $a>1$ and all $0\leq\lambda_1\leq 1.$

Define the function
$$f_{\lambda_1}(a)=\lambda_1^{s \alpha}+k^\alpha(1-\lambda_1)^{s \alpha}[(a-\lambda_1)/(1-\lambda_1)]^{(s/(1-s))\alpha }-k^\alpha a^{(s/(1-s))\alpha }.$$
Then the function is local fractional continuous on the
$(\lambda_1,\infty)$ and
$$g(\lambda_1)=f_{\lambda_1}(1)=\lambda_1^{s \alpha}+k^\alpha(1-\lambda_1)^{s \alpha}-k^\alpha.$$

It is easy to see that $g$ is  local fractional continuous on
$[0,1]$ and $g(1)=1^\alpha-k^\alpha<0^\alpha.$ So there is a number
$\lambda_{1_0},$  $0<\lambda_{1_0}<1,$ such that
$g(\lambda_{1_0})=f_{\lambda_{1_0}}(1)<0^\alpha.$ The local
fractional continuity of $f_{\lambda_{1_0}}$ yields that
$f_{\lambda_{1_0}}(a)<0^\alpha$ for a certain $a>1,$ i.e. inequality
(3.5) does not hold, which means that $f\not\in GK_s^2.$
\end{pf}

 \section{Conclusion}
 In the paper, we introduce the definitions of two kinds of generalized $s$-convex function on fractal sets and study the properties of these generalized $s$-convex functions. When $\alpha=1,$ these results  are the  classical situation.

\noindent{\bf Acknowledgments}\;\;{The authors would like to express
their gratitude to the reviewers for their very valuable comments.
And, this work is supported by the National Natural Science
Foundation of China (No. 11161042).}


\begin{thebibliography}{20}
\small
\bibitem{GRA}J. R. Jonathan and P. A. Mattew, ``Jensen's inequality predicts effects of environmental variation,"
 Trends in Ecology \& Evolution, vol. 14, no. 9, pp. 361-366, 1999.
\bibitem{GR} M. Grinalatt and J. T. Linnainmaa, ``Jensen's inequality, parameter uncertainty, and multiperiod investment,"
          Review of Asset Pricing Studies, vol. 1, no. 1, pp. 1-34, 2011.
\bibitem{HL} H. Hudzik and L. Maligranda, ``Some remarks on $s$-convex functions," Aequationes Mathematicae, vol. 48, no. 1, pp. 100-111, 1994.
\bibitem{S}E. Set, ``New inequalities of Ostrowski type for mappings whose derivatives are $s$-convex in the second sense via fractional
integrals," Computers \& Mathematics with Applications, vol.63, no.
7, pp. 1147-1154, 2012.
\bibitem{AK}M. Avci, H. Kavurmaci and M. E. Ozdemir, ``New inequalities of Hermite-Hadamard type via $s$-convex functions in the second sense with applications,"
 Applied Mathematics and Computation, vol. 217, pp. 5171-5176, 2011.
\bibitem{SYQ} S. Ye, H. Yin, and F. Qi,. ``Hermite-Hadamard type integral inequalities for geometric-arithmetically $s$-convex functions,"
 Analysis International Mathematical Journal of Analysis and Its Applications, vol. 33, no. 2, pp. 197-208, 2013.
\bibitem{KG} K. M. Kolwankar and A. D. Gangal, ``Local Fractional Calculus: a Calculus for Fractal Space-time," in Fractals: Theory and Application in Engineering, pp. 171-181, Springer, London, UK, 1999.
\bibitem{GB} A. K. Golmankhaneh and D. Baleanu, ``On a new measure on fractals," Journal of Inequalitiesand Applications,  vol. 522, no. 1, pp. 1-9, 2013.
\bibitem{M} B. B. Mandelbrot, The Fractal Geometry of Nature, Macmillan, New York, NY, USA, 1983.
\bibitem{F} K. J. Falconer, Fractal Geometry: Mathematical Foundations and Applications, Wiley, New York, NY, USA,  2007.
\bibitem{E} G. A. Edgar, Integral Probability and Fractal Measures, Springer, Berlin, GER, 1998.
\bibitem{KST}A. A. Kilbas, H. M. Srivastava and J. J. Trujillo, ``Theory and Applications of Fractional Differential Equations," North-Holland Mathematical Studies,
             vol. 204, Elsevier (North-Holland) Science Publishers, Amsterdam, London and New York, 2006.
\bibitem{Y} X. J. Yang, Advanced Local Fractional Calculus and Its Applications, World Science Publisher, New York, NY, USA, 2012.
\bibitem{MSY}H. Mo, X. Sui and D. Yu, ``Generalized convex functions on fractal sets and two related inequalities," Abstract and Applied Analysis, vol. 2014,
2014. Article ID 636751, 7 pages.
\bibitem{WSZW}W. Wie, H. M. Srivastava, Y. Zhang, L. Wang, P. Shen and J.Zhang, ``A local fractional integral inequality on fractal space
analogous to Anderson's inequality," Abstract and Applied Analysis,
vol. 2014, 2014. Article ID 797561, 7 pages.
\bibitem{YBM} X. J. Yang, D. Baleanu and J. A. T. Machado, ``Mathematical aspects of heisenberg uncertainty principle within local fractional fourier analysis,"
             Boundary Value Problems, vol. 2013, no. 1, pp. 131-146, 2013.
\bibitem{GYK}Y. Zhao, D. F. Cheng and X. J. Yang, ``Approximation solutions for local fractional Schr\"{o}dinger Equation in the on-dimensional cantorian system,"
  Advances in Mathematical Physics, vol. 2013, pp. 5, 2013.
\bibitem{YG1}A. M. Yang, Z. S. Chen, H. M. Srivastava and X. J. Yang, ``Application of the local fractional series expansion method and the variational iteration method
             to the Helmholtz equation involving local fractional derivative operators," Abstract and Applied Analysis, vol. 2013, 2013, Article ID
             259125, 6 pages.
\bibitem{YBKD} X. J. Yang, D. Baleanu, Y. Khan and  S. T. Mohyud-Din, ``Local fractional variational iteration method for diffusion and
              wave equations on Cantor sets," Romanian Journal of Physics, vol.59, no. 1-2, pp. 36-48, 2014.











\end{thebibliography}
\end{document}